# RECURRENCE OF RANDOM WALK TRACES[1]

By Itai Benjamini, Ori Gurel-Gurevich and Russell Lyons

*Weizmann Institute of Science, Weizmann Institute of Science and Indiana University*

We show that the edges crossed by a random walk in a network form a recurrent graph a.s. In fact, the same is true when those edges are weighted by the number of crossings.

**1. Introduction.** Let $G = (\mathsf{V}, \mathsf{E})$ be a locally finite graph and let $c: \mathsf{E} \to [0, \infty)$ be an assignment of *conductances* to the edges. We call $(G, c)$ a *network*. The associated random walk has transition probabilities $p(x, y) := c(x, y)/\pi(x)$, where $\pi(x) := \sum_y c(x, y)$. Assume that the network random walk is transient when it starts from some fixed vertex $o$. How big can the trace be, the set of edges traversed by the random walk? We show that they form a.s. a recurrent graph (for a simple random walk).

This fact is already known when $G$ is a Euclidean lattice and $c \equiv 1$ since a.s. the paths there have infinitely many cut-times, a time when the past of the path is disjoint from its future; see [7] and [8]. From this, recurrence follows by the criterion of Nash-Williams [12]. By contrast, Lyons and Peres [9] constructed an example of a transient birth-and-death chain which a.s. has only finitely many cut-times.

A result of similar spirit to ours was proved by Morris [11], who showed that the components of the wired uniform spanning forest are a.s. recurrent. For another a.s. recurrence theorem (for distributional limits of finite planar graphs), see [5].

We expect that a Brownian analogue of the theorem is true, that is, a.s. parabolicity of the Wiener sausage, with reflected boundary conditions. For background on recurrence in the Riemannian context, see, for example, [6]. It would be interesting to prove similar theorems for other processes. For example, consider the trace of a branching random walk on a graph $G$.

Received December 2005; revised June 2006.
[1]Supported in part by NSF Grant DMS-04-06017.

Then we conjecture that almost surely the trace is recurrent for a branching random walk with the same branching law. Perhaps a similar result holds for general tree indexed random walks. See Benjamini and Peres [3, 4] for definitions and background.

Perhaps one can strengthen our result as follows. Given a transient network $(G, c)$, denote by $T_n$ the trace of the first $n$ steps of the network random walk. Let $R(n)$ be the maximal effective resistance on $T_n$ between $o$ and another vertex of $T_n$, where each edge has unit conductance. By our theorem, $R(n) \uparrow \infty$ a.s. [Note, of course, that $R(n) \uparrow \infty$ for growing subgraphs does not imply recurrence of their union, as balls in the binary tree show.] Is there a uniform lower bound over all transient networks for the rate at which $R(n) \uparrow \infty$? That is, does there exist a function $f$ with $\lim_n f(n) = \infty$ such that for any transient network,

$$\limsup_n R(n)/f(n) > 0 \quad \text{a.s.?}$$

In particular, one can speculate that $f(n) = \log^2 n$ might work, which would arise from the graph $\mathbb{Z}^2$ (although $\mathbb{Z}^2$ is recurrent, it is on the border of transience). On the other hand, transient wedges in $\mathbb{Z}^3$ might allow one to prove that there is no such $f$.

**2. Proof.** Our proof will demonstrate the following stronger results. Let $N(x, y)$ denote the number of traversals of the edge $(x, y)$.

THEOREM 2.1. *The network $(G, \mathbf{E}[N])$ is recurrent. The networks $(G, N)$ and $(G, \mathbf{1}_{\{N>0\}})$ are a.s. recurrent.*

We shall use some facts relating electrical networks to random walks. See [10] for more background.

Let $\mathcal{G}(x, y)$ be the Green function, that is, the expected number of visits to $y$ for a network random walk started at $x$.

The *effective resistance* from a vertex $o$ to infinity is defined to be the minimum energy $\frac{1}{2} \sum_{x \neq y} \theta(x, y)^2 / c(x, y)$ of any unit flow $\theta$ from $o$ to infinity. This also equals

(2.1) $$\alpha := \mathcal{G}(o, o)/\pi(o).$$

In particular, the effective resistance is finite iff the network random walk is transient. Its reciprocal, *effective conductance*, is given by Dirichlet's principle as the infimum of the Dirichlet energy $\frac{1}{2} \sum_{x \neq y} c(x, y)[F(x) - F(y)]^2$ over all functions $F : \mathsf{V} \to [0, 1]$ that have finite support and satisfy $F(o) = 1$. Since the functional $c \mapsto \sum_{x \neq y} c(x, y)[F(x) - F(y)]^2$ is linear for any given $F$, we see that effective conductance is concave in $c$. Thus, if the conductances $\langle \mathbf{E}[N(x, y)]; (x, y) \in \mathsf{E}(G) \rangle$ give a recurrent network, then so a.s. do $\langle N(x, y);$

$(x, y) \in \mathsf{E}(G)\rangle$. Furthermore, Rayleigh's monotonicity principle implies that if $(G, N)$ is recurrent, then so is $(G, \mathbf{1}_{\{N>0\}})$. (Of course, it follows that any finite union of traces, whether independent or not, is also recurrent a.s.)

Thus, it remains to prove that $(G, \mathbf{E}[N])$ is recurrent. We shall, however, also show how the proof that $(G, N)$ is a.s. recurrent follows from a simpler argument. Another mostly probabilistic proof of this is due to Benjamini and Gurel-Gurevich [2].

The *effective resistance* from a finite set of vertices $A$ to infinity is defined to be the effective resistance from $a$ to infinity when $A$ is identified to a single vertex, $a$. The *effective resistance* from an infinite set of vertices $A$ to infinity is defined to be the infimum of the effective resistance from $B$ to infinity among all finite subsets $B \subset A$. Its reciprocal, the *effective conductance* from $A$ to infinity in the network $(G, c)$, will be denoted by $\mathcal{C}(A, G, c)$. From the above, we have

$$(2.2) \quad \mathcal{C}(A, G, c) = \sup_B \inf \left\{ \tfrac{1}{2} \sum_{x \neq y} c(x, y)[F(x) - F(y)]^2; \ F \restriction B \equiv 1, \ F \text{ has finite support} \right\},$$

where the supremum is over finite subsets $B$ of $A$.

Let the original voltage function be $v(\bullet)$ throughout this article, where $v(o) = 1$ and $v(\bullet)$ is 0 "at infinity." Then $v(x)$ is the probability of ever visiting $o$ for a random walk starting at $x$.

Note that

$$\mathbf{E}[N(x, y)] = \mathcal{G}(o, x)p(x, y) + \mathcal{G}(o, y)p(y, x)$$
$$= (\mathcal{G}(o, x)/\pi(x) + \mathcal{G}(o, y)/\pi(y))c(x, y)$$

and

$$\pi(o)\mathcal{G}(o, x) = \pi(x)\mathcal{G}(x, o) = \pi(x)v(x)\mathcal{G}(o, o).$$

Thus, we have [from the definition (2.1)]

$$(2.3) \quad \begin{aligned} \mathbf{E}[N(x, y)] &= \alpha c(x, y)[v(x) + v(y)] \\ &\leq 2\alpha \max\{v(x), v(y)\}c(x, y) \leq 2\alpha c(x, y). \end{aligned}$$

In a finite network $(H, c)$, we write $\mathcal{C}(A, z; H, c)$ for the effective conductance between a subset $A$ of its vertices and a vertex $z$. This is given by Dirichlet's principle as the infimum of the Dirichlet energy of $F$ over all functions $F : \mathsf{V}(H) \to [0, 1]$ that satisfy $F \restriction A \equiv 1$ and $F(z) = 0$. Clearly, $A \subset B \subset \mathsf{V}$ implies that $\mathcal{C}(A, z; H, c) \leq \mathcal{C}(B, z; H, c)$. The function that minimizes the Dirichlet energy is the voltage function, $v$. The amount of current that flows from $A$ to $z$ in this case is defined as $\sum_{x \in A, y \notin A}[v(x) - v(y)]c(x, y)$; it equals $\mathcal{C}(A, z; H, c)$. The voltage function that is $t$ on $A$ instead of 1 has Dirichlet energy equal to $t^2 \mathcal{C}(A, z; H, c)$ and gives a current that is $t$ times as large as the unit-voltage current, which shows that $\mathcal{C}(A, z; H, c)$ is the amount of current that flows from $A$ to $z$ divided by the voltage on $A$.

LEMMA 2.1. *Let $(H, c)$ be a finite network and $a, z \in \mathsf{V}(H)$. Let $v$ be the voltage function that is 1 at $a$ and 0 at $z$. For $0 < t < 1$, let $A_t$ be the set of vertices $x$ with $v(x) \geq t$. Then $\mathcal{C}(A_t, z; H, c) \leq \mathcal{C}(a, z; H, c)/t$. Thus, for every $A \subset \mathsf{V}(H) \setminus \{z\}$, we have*

$$\mathcal{C}(A, z; H, c) \leq \frac{\mathcal{C}(a, z; H, c)}{\min v \restriction A}.$$

PROOF. We subdivide edges as follows. If any edge $(x, y)$ is such that $v(x) > t$ and $v(y) < t$, then subdividing the edge $(x, y)$ with a vertex $z$ by giving resistances

$$r(x, z) := \frac{v(x) - t}{v(x) - v(y)} r(x, y)$$

and

$$r(z, y) := \frac{t - v(y)}{v(x) - v(y)} r(x, y)$$

will result in a network such that $v(z) = t$ while no other voltages change. Doing this for all such edges gives a possibly new graph $H'$ and a new set $A'_t$ whose internal vertex boundary is a set $W'_t$ on which the voltage is identically equal to $t$. We have $\mathcal{C}(A_t, z; H, c) = \mathcal{C}(A_t, z; H', c) \leq \mathcal{C}(A'_t, z; H', c)$. Now $\mathcal{C}(A'_t, z; H', c) = \mathcal{C}(a, z; H, c)/t$ since the amount of current that flows is $\mathcal{C}(a, z; H, c)$ and the voltage difference is $t$. Therefore, $\mathcal{C}(A_t, z; H, c) \leq \mathcal{C}(a, z; H, c)/t$, as desired.

For a general $A$, let $t := \min v \restriction A$. Since $A \subset A_t$, we have $\mathcal{C}(A, z; H, c) \leq \mathcal{C}(A_t, z; H, c)$. Combined with the inequality just reached, this yields the final conclusion. □

For $t \in (0, 1)$, let $\mathsf{V}_t := \{x \in \mathsf{V}; v(x) < t\}$. Let $W_t$ be the external vertex boundary of $\mathsf{V}_t$, that is, the set of vertices outside $\mathsf{V}_t$ that have a neighbor in $\mathsf{V}_t$. Write $G_t$ for the subgraph of $G$ induced by $\mathsf{V}_t \cup W_t$.

We will refer to the conductances $c$ as the *original* ones and the conductances $\mathbf{E}[N]$ as the *new* ones for convenience.

LEMMA 2.2. *The effective conductance from $W_t$ to $\infty$ in the network $(G_t, \mathbf{E}[N])$ is at most 2.*

PROOF. If any edge $(x,y)$ is such that $v(x) > t$ and $v(y) < t$, then subdividing the edge $(x,y)$ with a vertex $z$ as in the proof of Lemma 2.1 and consequently adding $z$ to $W_t$ has the effect of raising the conductance of the edge $(x,y)$ to $c(z,y) = c(x,y)[v(x) - v(y)]/[t - v(y)]$ and also, by (2.3), of raising its conductance in the new network from $\mathbf{E}[N(x,y)]$ to

$$\begin{aligned}\alpha c(z,y)[t + v(y)] &= \alpha c(z,y)[t - v(y) + 2v(y)] \\ &= \alpha c(x,y)[v(x) - v(y)] + 2\alpha c(z,y)v(y) \\ &> \alpha c(x,y)[v(x) - v(y)] + 2\alpha c(x,y)v(y) = \mathbf{E}[N(x,y)].\end{aligned}$$

Since raising edge conductances clearly raises effective conductance, it suffices to prove the lemma in the case that $v(x) = t$ for all $x \in W_t$. Thus, we assume this case for the remainder of the proof.

Suppose that $\langle (H_n, c); n \geq 1 \rangle$ is an increasing exhaustion of $(G,c)$ by finite networks that include $o$. Identify the boundary (in $G$) of $H_n$ to a single vertex, $z_n$. Let $v_n$ be the corresponding voltage functions with $v_n(o) = 1$ and $v_n(z_n) = 0$. Then $\mathcal{C}(o, z_n; H_n, c) \downarrow 1/\alpha$ and $v_n(x) \uparrow v(x)$ as $n \to \infty$ for all $x \in \mathsf{V}(G)$. Let $A$ be a finite subset of $W_t$. By Lemma 2.1, as soon as $A \subset \mathsf{V}(H_n)$, we have that the effective conductance from $A$ to $z_n$ of $H_n$ is at most $\mathcal{C}(o, z_n; H_n, c)/\min\{v_n(x); x \in A\}$. Therefore by Rayleigh's monotonicity principle, $\mathcal{C}(A, G_t, c) \leq \mathcal{C}(A, G, c) = \lim_{n \to \infty} \mathcal{C}(A, z_n; H_n, c) \leq 1/(\alpha t)$. Since this holds for all such $A$, we have

(2.4) $$\mathcal{C}(W_t, G_t, c) \leq 1/(\alpha t).$$

By (2.3), the new conductances on $G_t$ are obtained by multiplying the original conductances by factors that are at most $2\alpha t$. Combining this with (2.4), we obtain that the new effective conductance from $W_t$ to infinity in $G_t$ is at most 2. □

When the complement of $\mathsf{V}_t$ is finite for all $t$, which is the case for "most" networks, this completes the proof by the following lemma (and by the fact that $\bigcap_{t>0} \mathsf{V}_t = \varnothing$):

LEMMA 2.3. *If $H$ is a transient network, then for all $m > 0$, there exists a finite subset $K \subset \mathsf{V}(H)$ such that for all finite $K' \supseteq K$, the effective conductance from $K'$ to infinity is more than $m$.*

PROOF. Let $\theta$ be a unit flow of finite energy from a vertex $o$ to $\infty$. Since $\theta$ has finite energy, there is some $K \subset \mathsf{V}(G)$ such that the energy of $\theta$ on the edges with some endpoint not in $K$ is less than $1/m$. That is, the effective resistance from $K$ to infinity is less than $1/m$. □

Even when the complement of $\mathsf{V}_t$ is not finite for all $t$, this is enough to show that the network $(G, N)$ is a.s. recurrent: If $X_n$ denotes the position of the random walk on $(G, c)$ at time $n$, then $v(X_n) \to 0$ a.s. by Lévy's 0–1 law. Thus, the path is a.s. contained in $\mathsf{V}_t$ after some time, no matter the value of $t > 0$. By Lemma 2.3, if $(G, N)$ is transient with probability $p > 0$, then $\mathcal{C}(B_n, G, N)$ tends in probability, as $n \to \infty$, to a random variable that is infinite with probability $p$, where $B_n$ is the ball of radius $n$ about $o$. In particular, this effective conductance is at least $6/p$ with probability at least $p/2$ for all large $n$. Fix $n$ with this property. Let $t > 0$ be such that $\mathsf{V}_t \cap B_n = \varnothing$. Write $D$ for the (finite) set of vertices in $G$ incident to an edge $e \notin G_t$ with $N(e) > 0$. Then $\mathcal{C}(W_t, G_t, N) = \mathcal{C}(W_t \cup D, G, N) \geq \mathcal{C}(B_n, G, N)$. However, this implies that $\mathcal{C}(W_t, G_t, \mathbf{E}[N]) \geq \mathbf{E}[\mathcal{C}(W_t, G_t, N)] \geq 3$, which contradicts Lemma 2.2.

To complete the proof that $(G, \mathbf{E}[N])$ is recurrent in general, we show that although $\mathsf{V}_t$ may not separate the source $o$ from infinity, its complement in the network is recurrent:

LEMMA 2.4. *The vertices $\mathsf{V} \setminus \mathsf{V}_t$ induce a recurrent network for the original and for the new conductances.*

PROOF. Condition that the original random walk on $G$ returns to its starting point, $o$. Of course, the corresponding Doob-transformed Markov chain is recurrent. This corresponds to transformed transition probabilities $p(x, y)v(y)/v(x)$ for $x \neq o$, whence to transformed conductances $c'(x, y) := c(x, y)v(x)v(y)$. Rayleigh's monotonicity principle gives that when we delete $\mathsf{V}_t$, we still have a recurrent network. But off of $\mathsf{V}_t$, the conductances $c'$ differ by a bounded factor from the original conductances and also from the new conductances. This means that the part remaining after we delete $\mathsf{V}_t$ is recurrent for both the original and new conductances. □

PROOF OF THEOREM 2.1. The function $x \mapsto v(x)$ has finite Dirichlet energy for the original network, hence for the new (since conductances are multiplied by a bounded factor). Assume (for a contradiction) that the new random walk is transient. Then by Ancona, Lyons and Peres [1], $\langle v(X_n) \rangle$ converges a.s. for the new random walk. By Lemma 2.4, it a.s. cannot have a limit $> t$ for any $t > 0$, so it converges to 0 a.s.

This means that the unit current flow $i$ for the new network (which is the expected number of signed crossings of edges) has total flow 1 through $W_t$ into $G_t$ for all $t > 0$. Thus, we may choose a finite subset $A_t$ of $W_t$ through which at least $1/2$ of the new current enters. With the notation $(d_t^* i)(x) := \sum_{y \in \mathsf{V}(G_t)} i(x, y)$, this means that $\sum_{x \in A_t} d_t^* i(x) \geq 1/2$. By Lemma 2.2, there is a function $F_t : \mathsf{V}_t \cup W_t \to [0, 1]$ with finite support and with $F_t \equiv 1$ on $A_t$ whose Dirichlet energy on the network $(G_t, \mathbf{E}[N])$ is at most 3. Write

$(dF_t)(x,y) := F_t(x) - F_t(y)$. By the Cauchy–Schwarz inequality, we have

$$\left[\sum_{x \neq y \in \mathsf{V}(G_t)} i(x,y)\, dF_t(x,y)\right]^2 \leq \sum_{x \neq y \in \mathsf{V}(G_t)} i(x,y)^2/c(x,y)$$
$$\times \sum_{x \neq y \in \mathsf{V}(G_t)} c(x,y)\, dF_t(x,y)^2$$
$$\leq 3 \sum_{x \neq y \in \mathsf{V}(G_t)} i(x,y)^2/c(x,y).$$

On the other hand, summation by parts yields that

$$\sum_{x \neq y \in \mathsf{V}(G_t)} i(x,y)\, dF_t(x,y) = \sum_{x \in \mathsf{V}(G_t)} d_t^* i(x) F_t(x) \geq \sum_{x \in A_t} d_t^* i(x) \geq 1/2.$$

Therefore, $\sum_{x \neq y \in \mathsf{V}(G_t)} i(x,y)^2/c(x,y) \geq 1/12$, which contradicts $\bigcap_t \mathsf{V}(G_t) = \varnothing$ and the fact that $i$ has finite energy. $\square$

**Acknowledgments.** Thanks to Gidi Amir, Gady Kozma, Ron Peled and Benjy Weiss for useful discussions. We also thank the referees for useful comments.

I. BENJAMINI
O. GUREL-GUREVICH
MATHEMATICS DEPARTMENT
THE WEIZMANN INSTITUTE OF SCIENCE
REHOVOT 76100
ISRAEL
E-MAIL: itai@wisdom.weizmann.ac.il
   origurel@weizmann.ac.il
URL: http://www.wisdom.weizmann.ac.il/~itai/
   http://www.wisdom.weizmann.ac.il/~origurel/

R. LYONS
DEPARTMENT OF MATHEMATICS
INDIANA UNIVERSITY
BLOOMINGTON, INDIANA 47405
USA
E-MAIL: rdlyons@indiana.edu
URL: http://mypage.iu.edu/~rdlyons/